\documentclass[12pt]{article}

\usepackage[T1]{fontenc}
\usepackage{csquotes,enumerate,amssymb,amsmath}
\usepackage[backend=bibtex,firstinits=true, url=false, isbn=false]{biblatex}
\usepackage{amsmath}
\def\theequation{\arabic{equation}}

  \def\theequation{\thesection.\arabic{equation}}
  \def\theequation{\arabic{equation}}

\newsavebox{\savepar}

\newtheorem{remark}{Remark}[section]

\def\appendix{\par
   \setcounter{section}{0}
   \setcounter{equation}{0}
   \def\@chapapp{APPENDIX}
   \def\thesection{\Alph{section}}
   \def\theequation{A.\arabic{equation}}
   \def\section{
   \setcounter{equation}{0}
   \refstepcounter{section}
   \@startsection {section}{A}{0pt}{-3.5ex \@plus -1ex \@minus -.2ex}
                  {0.3ex \@plus.2ex}{\normalsize\sffamily\bfseries} }}

\def\eop{\hbox{{\vrule height7pt width3pt depth0pt}}}

\makeatletter
\newcommand{\least}{\let\CS=\@currsize\renewcommand{\baselinestretch}{1}\tiny\CS}
\newcommand{\singlespacing}{\renewcommand{\baselinestretch}{1.5}}
\newcommand{\oneandahalfspacing}{\let\CS=\@currsize\renewcommand{\baselinestretch}{1.2}\tiny\CS}
\newcommand{\doublespacing}{\let\CS=\@currsize\renewcommand{\baselinestretch}{2.5}\tiny\CS}

\usepackage{amsmath}
\usepackage{amsfonts}
\usepackage{euscript}
\usepackage{oldgerm}
\usepackage{rotating}

  \renewcommand{\baselinestretch}{1.3}

\makeindex

\numberwithin{equation}{section}
\newcommand{\be}{\begin{equation}}
\newcommand{\ee}{\end{equation}}

\newcommand{\beano}{\begin{eqnarray*}}
\newcommand{\eeano}{\end{eqnarray*}}
\newcommand{\bea}{\begin{eqnarray}}
\newcommand{\eea}{\end{eqnarray}}

\newcommand{\ba}{\begin{array}}
\newcommand{\ea}{\end{array}}

\newcommand{\vone}{\vskip 2ex}
\newcommand{\vtwo}{\vskip 4ex}

\begin{document}
\singlespacing

\begin{center}
{\Large \bf Exponential Order Statistics, the Basel problem  and  Combinatorial Identities}\\
\end{center}

\vtwo
\begin{center}
{\large \bf  P. Vellaisamy\textsuperscript{1} and A. Zeleke \textsuperscript{2}} \\
 \textsuperscript{1}Department of Mathematics, Indian Institute of Technology Bombay \\
\noindent Powai, Mumbai-400076, India.\\
{Email: \it pv@math.iitb.ac.in}

 \textsuperscript{2}Lyman Briggs College \& Department of Statistics \& Probability, \\
\noindent Michigan State University. East Lansing, MI 48825, USA\\
{Email: \it zeleke@stt.msu.edu}
\end{center}

\abstract {We consider the $k$-th order statistic from unit exponential distribution and show that it can be represented as a sum of
independent exponential random variables. Our proof is simple and different. It readily proves that the standardized exponential spacings also follow unit exponential distribution. An interesting probabilistic proof of the Basel problem is also given. Another advantage of our approach is that by computing
the Laplace transform of the $k$-th order statistic in two different ways, we derive  several interesting combinatorial identities. A probabilistic interpretation of these identities and their generalizations are also   given.}

\vone

\vone
\noindent  MSC2010 {\it Subject Classification}: Primary: 62E15, 05A19; Secondary: 60C05.\\

\section{Introduction}
It is known that order statistics from exponential distribution have several interesting properties. We consider, without loss of generality, the exponential distribution with mean unity. For example, the $k$-th order statistic has the distribution of sum of independent exponential random variables (rvs). Another interesting result is that the spacings of order statistics also follow exponential distribution. The usual proofs of these results use the transformation to the set of spacings from the set of order statistics and by applying Jacobian density theorem.

In this paper, we prove the above-mentioned and some other
results using the Laplace transform methods. This approach is simpler and indeed several related results can be proved in a unified way. We provide also an interesting probabilistic proof of the Basel problem.
Another key purpose of this article is to bring out the connection between exponential order statistics and several combinatorial identities.
In fact, we give simpler proofs of several combinatorial/binomial identities by evaluating  the Laplace transformation of the $k$-th exponential order statistic by two different ways and equating them. We also point out the probabilistic  interpretations of these combinatorial identities. 

\section{The Order Statistics Results} Let $X_{1},\ldots, X_{n}$ be iid continuous rvs with cdf $F(x)$ and density $f(x)$.
Also, let $X_{(1)}<X_{(2)}<\ldots X_{(n)}$ be the corresponding order statistics. Then it is well known (see David and Nagaraja (2003)) that the distribution and the density of $k$-th order statistic $X_{(k)}$ are respectively given by
\begin{equation}\label{order2}
F_{k}(x)=P(X_{(k)}<x)=\sum_{m=k}^{n}\binom{n}{m} F^{m}(x)(1-F(x))^{n-m},
\end{equation}
and
\begin{equation}\label{order3}
f_{k}(x)=n\binom{n-1}{k-1} F^{k-1}(x)(1-F(x))^{n-k} f(x),
\end{equation}
for $1\leq k\leq n$.

 \noindent Our interest is on the exponential order statistics. Let $T_{1},\ldots,T_{n}$ be independent unit exponential  rvs
 so that
 \begin{equation}\label{order1}
 	f_{T_{1}}(t)=e^{-t}, ~~t>0.
 \end{equation}
 and let $T_{(1)}<T_{(2)}<\ldots T_{(n)}$ be the order statistics of $T_{j}$'s.
Then the density of $T_{(k)}$ is
\begin{equation}\label{order4}
f_{k}(t)=\frac{1}{B(k,n-k+1)}\left(1-e^{-t}\right)^{k-1}e^{-(n-k+1)t},~~t>0.
\end{equation}

\noindent Here $B(m,n)$ is the usual beta function.

First we prove an interesting property of $T_{(k)}$ by computing the Laplace transform of $T_{(k)}$.

\noindent Note for $s>0$,
\begin{equation*}
\mathbb{E}(e^{-sT_{(k)}})=\frac{1}{B(k,n-k+1)}\int_{0}^{\infty}e^{-(s+1)t}\left(1-e^{-t}\right)^{k-1} e^{-(n-k)t}dt.
\end{equation*}
Substitute $w=\left(1-e^{-t}\right)$ to obtain
\begin{align}
\mathbb{E}(e^{-sT_{(k)}})&=\frac{1}{B(k,n-k+1)}\int_{0}^{1} w^{k-1} (1-w)^{n-k+s+1-1} dw \nonumber \\
&=\frac{B(k,n-k+s+1)}{B(k,n-k+1)}\nonumber \\
&=\frac{\Gamma(n+1)\Gamma(n-k+s+1)}{\Gamma(n-k+1)\Gamma(n+s+1)}\nonumber \\
&=\frac{n(n-1)(n-2)\ldots (n-k+1)}{(s+n)(s+n-1)\ldots (s+n-k+1)}\nonumber\\
&=\prod_{j=n-k+1}^{n}\left(\frac{j}{s+j}\right)\nonumber\\
& =f_{n, k}(s)~~\text{(say)}.\label{order5}
\end{align}
Note also that
\begin{align}
\mathbb{E}(e^{-sT_{(k)}})&=\prod_{j=n-k+1}^{n}\mathbb{E}(e^{-sY_{j}})\nonumber \\
&=\prod_{j=n-k+1}^{n}\mathbb{E}(e^{-s\sum_{j=n-k+1}^{n}Y_{j}}),\label{order6}
\end{align}
where $Y_{1},\ldots, Y_{n}$ are independent exponential rvs and $Y_{j}\sim \text{Exp}(j)$ with mean $j^{-1}$, $1\leq j\leq n$.
Thus,  for $1\leq k\leq n$, we have from \eqref{order6}
\begin{equation}\label{order7}
T_{(k)}\stackrel{d}{=}\sum_{j=n-k+1}^{n}Y_{j},
\end{equation}
where $X\stackrel{d}{=}Y$ means both $X$ and $Y$ have identical distributions.

\noindent Note when $k=n$, we have
\begin{equation*}
T_{(n)}\stackrel{d}{=}\sum_{j=1}^{n}Y_{j},
\end{equation*}
a known result (see, for example, (4.1) of Vellaisamy (2015)).

\noindent Also, when $k=1$,
\begin{equation*}
T_{(1)}\stackrel{d}{=}Y_{n},
\end{equation*}
where $Y_{n}\sim \text{Exp}(n)$, a well-known result.

\noindent Next some interesting remarks are in order.
\begin{remark} {\em
	From (\ref{order7}), we obtain a known result that
	\begin{align}
	T_{(k)}&\stackrel{d}{=}Y_{n-k+1}+Y_{n-k+2}+\ldots+Y_{n}\nonumber\\
	&\stackrel{d}{=} W_{1}+W_{2}+\ldots+W_{k},\label{order8}
	\end{align}
	where $W_{j}$'s are independent and $W_{j}\sim \text{Exp}(n-k+j)$.
	Also, from (\ref{order8}), we get
	\begin{equation}\label{order9}
	\mathbb{E}(T_{(k)})=\sum_{j=1}^{k}\mathbb{E}(W_{j})=\sum_{j=1}^{k}\frac{1}{(n-k+j)}
	\end{equation}
	and
	\begin{equation}\label{order10}
	\mathbb{V}ar(T_{(k)})=\sum_{j=1}^{k}\mathbb{V}ar(W_{j})=\sum_{j=1}^{k}\frac{1}{(n-k+j)^2},
	\end{equation}
	for $1\leq k\leq n$.}
\end{remark}
Note that usual proofs of the  results in \eqref{order9} and
\eqref{order10}, based on the density of $T_{(k)}$, are rather complicated.

\begin{remark} {\em
	Also, from the representation given in (\ref{order7}), we can immediately obtain the distribution of the spacings $(T_{(k)}-T_{(k-1)})$, $1\leq k\leq n$, with $T_{(0)}\equiv 0$.

	\noindent Since
	\begin{equation*}
	T_{(k)}\stackrel{d}{=}Y_{n-k+1}+Y_{n-k+2}+\ldots+Y_{n}
	\end{equation*}
	we have for $1\leq k\leq n $,
	\begin{align*}
	T_{(k)}-T_{(k-1)}&\stackrel{d}{=}Y_{n-k+1}
		\stackrel{d}{=} \frac{T_{1}}{(n-k+1)}
	\end{align*}
	or equivalently the normalized spacings
	\begin{equation}\label{order11}
	(n-k+1)(T_{(k)}-T_{(k-1)})\stackrel{d}{=} T_{1},~~1\leq k\leq n,
	\end{equation}
where $T_{1}\sim \text{Exp}(1)$. }
\end{remark}
\noindent The proof of the above result is usually based on the joint density of $T_{(1)},\ldots, T_{(k)}$ and then making suitable transformations to spacings. This is rather involved, whereas our proof easily follows from (\ref{order7}).

 \subsection{The Basel Problem}  In this section, we give an interesting probabilistic proof of the Euler's constant $\gamma$ and  the Basel problem.
	  Let $Z_n= (T_{(n)} - \ln(n)) $. Since,
\begin{equation*}
 P(Z_n \leq x) = (1-e^{-(\ln(n)+x)})^n
                 = (1- \frac{e^{-x}}{n})^n \longrightarrow \exp\{-e^{-x}\}, ~\text{as}~ n \rightarrow \infty,
\end{equation*}		
we have $Z_n  \overset {\cal L} {\Longrightarrow} Z $ (see Billingsley (1995, p.~329) for the definition of convergence in distribution),
where $Z$ follows standard Gumbel distribution with density

\begin{equation*}
 f(z)= e^{-z}\exp{(-e^{-z})},~  z \in \mathbb{R}.
\end{equation*}
It is known that $E(Z)= \nu$ (Euler's constant) and  $Var(Z)= \frac{\pi^2}{6}$.

\noindent Let now $G_n$ denote the cdf of $Z_n$. Then it can be shown that
\begin{eqnarray*}
	G_n(-x)+1-G_n(x) &=&
	\left\{
	\begin{array}{ll}
		(1-\frac{e^x}{n})^n + 1- (1-\frac{e^{-x}}{n})^n, & ~ \mbox{if}~ x < \log(n)\\
		1- (1-\frac{e^{-x}}{n})^n, & ~ \mbox {if}~  x \geq \log(n)\\
	\end{array}
	\right. \\
	&<&
	\left\{
	\begin{array}{ll}
		exp(-e^{x}) + 1- (1-e^{-x}), & ~ \mbox{if}~ x < \log(n)\\
		1- (1-e^{-x}), & ~ \mbox {if}~  x \geq \log(n),\\
	\end{array}
	\right.
\end{eqnarray*}
since h(x)= $(1- \frac{\theta}{x})^x$ increases to $e^{-\theta}$ for $\theta < x< \infty.$

\noindent Hence, for all $n \geq 1$, and for all $x>0$,
\begin{eqnarray} \label{neqn2.12}
	\sup_{n} \{G_n(-x)+1-G_n(x)\} < \max\{exp(-e^{x}) + 1-e^{-x}, e^{-x} \} < 2e^{-x}.	
\end{eqnarray}

\noindent Note that (see Pakshirajan (2003)p.~119) $|Z_n|^{k}$ is uniformy integrable if and only if
\begin{eqnarray*}
	\sup_{n} \int_{R}^{\infty} x^{k-1} [G_n(-x)+1-G_n(x)]dx
	\rightarrow 0,
\end{eqnarray*}
as $R \rightarrow \infty.$

\noindent Using \eqref{neqn2.12}, we get
\begin{eqnarray*}
	\sup_{n} \int_{R}^{\infty} x^{k-1} [G_n(-x)+1-G_n(x)]dx
	&<& 2\int_{R}^{\infty} x^{k-1} e^{-x} dx \longrightarrow 0,
\end{eqnarray*}
as $R \rightarrow \infty$.

\noindent Hence, $|Z_n|^{r}$ is uniformy integrable and hence
for  $r \geq 1 $,
$$\lim_{n\rightarrow \infty} E(Z_n^r) = E(Z^r).$$
Using the first and second moments convergence, we have
\begin{equation}
 \lim_{n\rightarrow \infty} E(Z_n)=\lim_{n\rightarrow \infty} (\sum_{j=1}^{n} \frac{1}{j} - \ln(n))= \gamma.
\label{neqn3.1}
\end{equation}
and
\begin{equation}
\lim_{n\rightarrow \infty} Var(Z_n)=\lim_{n\rightarrow \infty} (\sum_{j=1}^{n} \frac{1}{j^2})= \sum_{j=1}^{\infty} \frac{1}{j^2}= \frac{\pi^2}{6},
\label{neqn3.2}
\end{equation}
which is the famous Basel problem.
\noindent Note that other mathematical proofs for the equality of the last two terms of \eqref{neqn3.1} and
\eqref{neqn3.2} are rather involved. More importantly, it is interesting to observe the connections between the  exponential order statistics, Euler constant and the Basel problem. For other probabilistic proofs of \eqref{neqn3.2}, see Pace (2011) and Holst (2013).	

\subsection{Combinatorial Identities} In this section, we derive several combinatorial identities using the Laplace transformation of the $k$-th order statistic $T_{(k)}$. We derive  the Laplace transform of $T_{(k)}$ in two different ways and equate them to get identities. Note that we have already derived one form, using the pdf of $T_{(k)}$. We  obtain next another form using the cdf of $T_{(k)}$ and integration by parts.

\noindent First note that the cdf of $T_{(k)}$ from (\ref{order2}) is
\begin{equation}\label{order12}
F_{k}(x)=\sum_{m=k}^{n}\binom{n}{m} \left(1-e^{-t}\right)^{m}e^{-(n-m)t},~~t>0.
\end{equation}
Using (\ref{order12}), we have for $1\leq k\leq n$,
\begin{align}
\mathbb{E}(e^{-sT_{(k)}})&=\int_{0}^{\infty}e^{-st}dF_{k}(t)\nonumber\\
&=s\int_{0}^{\infty}e^{-st}F_{k}(t)dt~~~~~\text{(integration by parts)}\nonumber\\
&=s\int_{0}^{\infty}e^{-st}\left\{\sum_{m=k}^{n}\binom{n}{m} \left(1-e^{-t}\right)^{m}e^{-(n-m)t}\right\}dt\nonumber\\
&=s\sum_{m=k}^{n}\binom{n}{m}\int_{0}^{\infty}e^{-st}\left\{\sum_{j=0}^{m}\binom{m}{j}(-1)^{j} e^{-jt}\right\}e^{-(n-m)t}dt\nonumber\\
&=s\sum_{m=k}^{n}\binom{n}{m}\sum_{j=0}^{m}(-1)^{j}\binom{m}{j} \int_{0}^{\infty}e^{-(s+j+n-m)t}dt\nonumber \\
&=\sum_{m=k}^{n}\sum_{j=0}^{m}(-1)^{j}\binom{n}{m}\binom{m}{j}\left(\frac{s}{s+j+(n-m)}\right)\nonumber\\
&=g_{n, k}(s),~~(\text{say}).\label{order13}
\end{align}
Thus, from (\ref{order5}) and (\ref{order13}), we obtain an interesting combinatorial identity
\begin{eqnarray}\label{orderbinom}
\sum_{m=k}^{n}\sum_{j=0}^{m}(-1)^{j}\binom{n}{m}\binom{m}{j}\left(\frac{s}{s+n-m+j}\right) &=& \prod_{j=n-  
k+1}^{n}\left(\frac{j}{s+j}\right)
\end{eqnarray}
for  $1\leq k\leq n$, and $s>0$.

\noindent When $k=1$, we get from (\ref{orderbinom}),
\begin{equation}\label{orderbinom1}
\sum_{m=1}^{n}\sum_{j=0}^{m}(-1)^{j}\binom{n}{m}\binom{m}{j}\left(\frac{s}{s+n-m+j}\right)=\frac{n}{s+n},
\end{equation}
for $n\geq 1$ and $s>0$, an interesting binomial identity.

\noindent Similarly when $k=n$, we obtain from (\ref{orderbinom}),
\begin{eqnarray}\label{orderbinom2}
\sum_{j=0}^{n}(-1)^{j}\binom{n}{j}\left(\frac{s}{s+j}\right)=\prod_{j=1}^{n} \left(\frac{j}{s+j}\right). 
\end{eqnarray}
 
\begin{remark} {\em
 
Note that, when $s=k$ is a positive integer, (2.19) leads to $$\sum_{j=0}^{n}(-1)^{j}\binom{n}{j}\left(\frac{k}{k+j}\right)=\frac{1}{\binom{n+k} {k}}.$$

}
\end{remark}

\begin{remark} {\em
\noindent Recently,  Peterson (2013) derived the above binomial identity (2.19) from probabilistic considerations.  Vellaisamy (2015) gave a simple and different proof, based on Laplace transforms. He also discussed a statistical application of the above identity, using the binomial inversion formula (see Aigner (2007),~p. 73). Note also that
an application of binomial inversion to \eqref{orderbinom2}
leads to equation (4.7) of Vellaisamy (2015). Similarly, the equation
\eqref{orderbinom} will lead to a new identity.
}
\end{remark}

\begin{remark}
	Note from (\ref{orderbinom}),
	\begin{equation*}
	\sum_{m=k}^{n}\binom{n}{m}\left(\frac{s}{s+n-m}\right)\left\{\sum_{j=0}^{m}(-1)^{j}\binom{m}{j}\left(\frac{s+n-m}{s+n-m+j}\right)\right\}=\prod_{j=n-k+1}^{n}\left(\frac{j}{s+j}\right),
	\end{equation*}
	which using (\ref{orderbinom2}), leads to
	\begin{equation*}
	\sum_{m=k}^{n}\binom{n}{m}\left(\frac{s}{s+n-m}\right)\prod_{i=1}^{m}\left(\frac{i}{s+n-m+i}\right)=\prod_{j=n-k+1}^{n}\left(\frac{j}{s+j}\right),
	\end{equation*}
	for $1\leq k\leq n$, and $s>0$.
\end{remark}
\section{Probabilistic Connections and Generalizations} In this section, we give a probabilistic interpretation to the main combinatorial identity
$f_{n,k}(s)=g_{n,k}(s)$
for $n\geq 1$ and $s>0$.

\vspace*{-0.3cm}
\subsection{Probabilistic Interpretations}
\noindent Let now $X_{r}\sim G(s,r)$, $s>0$, $r\geq 1$ be a positive integer, with density
\begin{equation}\label{order14}
f(x|s,r)=\frac{s^r}{\Gamma(r)} e^{-sx} x^{r-1},~~x>0.
\end{equation}
 Then, it is well-known (see eq. (3.3.9) of Casella and Berger (2002)) that
\begin{equation}\label{order15}
\mathbb{P}(X_{r}>x)=\sum_{j=0}^{r-1}\frac{e^{-sx}(sx)^{j}}{j!}.
\end{equation}
Assume that $X_r$ and $T_{(k)}$ are independent.
Using \eqref{order15}, we obtain
\begin{align}
\mathbb{P}(X_{r}>T_{(k)})&=\mathbb{E}\left[\mathbb{P}(X_{r}>T_{(k)})|T_{(k)}\right]\nonumber\\
&=\sum_{j=0}^{r-1}\frac{s^{j}}{j!}\mathbb{E}\left(e^{-sT_{(k)}}T_{(k)}^{j}\right). \label{order16}
\end{align}

\noindent Note that $f_{n,k}(s)=\mathbb{E}(e^{-sT_{(k)}})$ satisfies
\begin{equation}\label{order17}
f^{(j)}_{n,k}(s)=(-1)^{j}\mathbb{E}\left(e^{-sT_{(k)}}T_{(k)}^{j}\right),
\end{equation}
where $h^{(j)}(s)$ denotes the $j$-th derivative of $h(s)=h^{(0)}(s)$.

\noindent Thus,  for $r\geq 1$ and $s>0$, we obtain from \eqref{order16} and \eqref{order17},
\begin{align}
\mathbb{P}(X_{r}>T_{(k)})&=\sum_{j=0}^{r-1}(-1)^{j}\frac{s^{j}}{j!}f^{(j)}_{n,k}(s)\label{order18}\\
&=\sum_{j=0}^{r-1}(-1)^{j}\frac{s^{j}}{j!}g^{(k)}_{n,k}(s), \label{order19}
\end{align}
since $f_{n,k}(s)=g_{n,k}(s)$.

\noindent When $r=1$, we obtain
\begin{equation}\label{order20}
f_{n,k}(s)=g_{n,k}(s)=P(X_{1}>T_{(k)}),	
\end{equation}
which shows that the binomial identity (\ref{orderbinom}) admits a probabilistic interpretation.

\noindent Similarly, when $r=2$, we get
\vspace*{-0.4cm}
\begin{equation}\label{order21}
f_{n,k}(s)-sf^{'}_{n,k}(s)=g_{n,k}(s)-sg^{'}_{n,k}(s)
\end{equation}
which is a new identity (see \eqref{eqn413}) and this corresponds to $P(X_{2}>T_{k})$.

\noindent Note however the combinatorial identity obtained from $f^{'}_{n,k}(s)=g^{'}_{n,k}(s)$ may not represent a probability of an event (see Vellaisamy (2015), p.~243).

\vspace*{-.3cm}

\subsection {Generalizations}

In this section, we generalize some of the binomial identities derived in Section 3.
Assume, as before, that $X_r$ and $T_{(k)}$ are independent. Observe that we computed, in the previous section, $\mathbb{P}(T_{(k)} <X_{r})$ by conditioning on $T_{(k)}$. We next compute it by conditioning on $X_r$.
Note that
\begin{align}
\mathbb{P}(T_{(k)} <X_{r})&=\mathbb{E}\left[\mathbb{P}(T_{(k)} < X_{r})|X_{r}\right]\nonumber\\
&=\mathbb{E}\left(\sum_{m=k}^{n} \binom{n}{m}(1-e^{-X_{r}})^{m} e^{-(n-m)X_{r}} \right) ~ (\text{from}~ \eqref{order12}) \nonumber \\
&=\sum_{m=k}^{n} \binom{n}{m} \left(\sum_{j=0}^{m} (-1)^{j}\binom{m}{j}\mathbb{E} e^{-(n-m+j)X_{r}} \right) \label{eqn16}	
\end{align}
Since $X_r \sim G(s, r)$, we have

\vspace*{-0.4cm}
\begin{align}
\mathbb{E} e^{-(n-m+j)X_{r}} &= \left( \frac{s}{s+n-m+j}\right)^{r}.
\end{align}

\noindent Hence, we get
\begin{align}
\mathbb{P}(T_{(k)} <X_{r})&=\sum_{m=k}^{n} \binom{n}{m} \sum_{j=0}^{m} (-1)^{j}\binom{m}{j}\left( \frac{s}{s+n-m+j}\right)^{r}. \label{eqn411}
\end{align}

\noindent Thus, we obtain from \eqref{order18} and \eqref{eqn411},
\begin{align} \label{eqn412}
\sum_{m=k}^{n} \binom{n}{m} \sum_{j=0}^{m} (-1)^{j}\binom{m}{j}\left( \frac{s}{s+n-m+j}\right)^{r} &=\sum_{j=0}^{r-1}(-1)^{j}\frac{s^{j}}{j!}f^{(j)}_{n,k}(s)
\end{align}
where,  as seen before,
\begin{equation*}
f_{n,k}(s)=\prod_{j=n-k+1}^{n}\left(\frac{j}{s+j}\right).
\end{equation*}

\noindent When $r=1$,  the identity in \eqref{eqn412}  reduces to \eqref{orderbinom}.

\noindent When  $r=2$, we get
\begin{align}
\sum_{m=k}^{n} \binom{n}{m} \sum_{j=0}^{m} (-1)^{j}\binom{m}{j}\left( \frac{s}{s+n-m+j}\right)^{2} &=f_{n,k}(s) -s f^{(1)}_{n,k}(s)\nonumber \\
&=\prod_{j=n-k+1}^{n}\left(\frac{j}{s+j}\right)\left[ 1+ \sum_{j=n-k+1}^{n}\frac{s}{s+j} \right] \label{eqn413}
\end{align}

\noindent When $k=n$, the above identity reduces to equation (4.3) of Vellaisamy (2015). However, when $k=1$,  the above result gives us a new identity, namely,
\begin{align}
\sum_{m=1}^{n} \binom{n}{m} \sum_{j=0}^{m} (-1)^{j}\binom{m}{j}\left( \frac{s}{s+n-m+j}\right)^{2} &= \frac{n(n+2s)}{(s+n)^2},
 \label{eqn414}
\end{align}
for $s>0$ and $n \geq 1$.
 Thus, we have generalized the basic binomial identity in \eqref{orderbinom2} in several directions and hope these  could be useful to other applied areas as well.

\vspace*{.2cm}
\noindent {\bf Acknowledgements}. This work was completed while the first author was visiting the department of statistics and probability, Michigan State University during Summer-2017. Also, the authors are grateful to Professor M. Sreehari for his comments and especially  in proving the uniform integrability result.

\vone
\noindent{\bf{\Large References}} 

\noindent Aigner, M. (2007). {\it A Course in Enumeration}. Berlin, Springer-Verlag.

\noindent Billingsley, P. (1995).  {\it Probability and Measure}. Third edition, John Wiley \& Sons, New York.

\noindent David, H. A. and Nagaraja, H. N. (2003). {\it Order Statistics}. Third edition, Wiley, New York.

\noindent Casella, G. and Berger, R. L. (2002). {\it Statistical Inference}. Second edition,  the Wadsworth \& Brooks, CA.

\noindent Holst, L. (2013). Probabilistic proofs the Euler identities. {\it Journal of Applied Probability}, {\it 50}, 1206-1212.

\noindent Pace, L. (2011). Probabilistically proving $\zeta(2)=
\frac{\pi^2}{6}$. {\it American Mathematical Monthly}, 118, 641-643.

\noindent Peterson, J. (2013). A probabilistic proof of a binomial identity. {\it American Mathematical Monthly}, 120, 558-562.

\noindent Pakshirajan, R. P. (2013). {\it Probability Theory:A foundational Course}. Hindustan Book Agency, New Delhi, India.

\noindent Vellaisamy, P. (2015).  On probabilistic proofs of certain binomial identities. {\it American Statistician}, 69, 241-243.

\vone
\end{document}